\documentclass[12pt,twoside]{article}
\usepackage{amsmath}
\usepackage{amssymb}
\usepackage[]{graphicx}
\usepackage{epsfig, float}
\usepackage[cp1251]{inputenc}
\usepackage[english]{babel}
\usepackage{color}

\newcommand{\be}{\begin{equation}}
\newcommand{\ee}{\end{equation}}
\newcommand{\bea}{\begin{eqnarray}}
\newcommand{\eea}{\end{eqnarray}}

\newtheorem{theorem}{Theorem}
\newtheorem{remark}{Remark}

\newtheorem{lemma}{Lemma}

\begin{document}

 \date{}

\title{Face reduction and the immobile indices approaches  to regularization of  linear Copositive Programming problems}

\author{Kostyukova O.I.\thanks{Institute of Mathematics, National Academy of Sciences of Belarus, Surganov str. 11, 220072, Minsk,
 Belarus  ({\tt kostyukova@im.bas-net.by}).} \and Tchemisova T.V.\thanks{Mathematical Department, University of Aveiro, Campus
Universitario Santiago, 3810-193, Aveiro, Portugal ({\tt tatiana@ua.pt}).}}
\maketitle

\begin{abstract}

 The paper is devoted to the regularization of linear Copositive Programming problems which consists  of  transforming a problem to an equivalent form, where  the Slater condition is satisfied and  the strong duality holds. We  describe here two regularization algorithms   based on the concept of immobile indices and an understanding of the important role  these indices play in  the feasible sets' characterization.   These   algorithms are compared to some    regularization procedures  developed for a more general case of convex problems and  based on a  facial reduction approach. We show that the immobile-index-based  approach  combined with   the specifics of copositive problems allows  us  to construct more explicit and detailed regularization algorithms for  linear Copositive Programming  problems  than  those   already available.

\end{abstract}

\textbf{Key words.} Linear copositive programming, strong duality, normalized immobile index set, regularization, minimal cone,  facial reduction, constraint qualifications
\\

\textbf{AMS subject classification.} 90C25, 90C30, 90C34

\section{Introduction}\label{Introduction} Conic optimization is a subfield of convex optimization that studies the problems of minimizing a convex function over the intersection of an affine subspace and a convex cone.  For a gentle introduction to conic optimization and  a survey of  its  applications in Operations Research and related areas, we refer interested readers to \cite{Let} and the references therein.

Copositive Programming (CoP) problems form a special class of conic  problems and can be considered as an optimization over the convex cone of so-called {\it copositive matrices} (i.e. matrices which are positive semi-defined on the non-negative orthant). Copositive models arise in many important applications, including $\mathcal{NP}$ -hard problems. For the references on motivation and application  of CoP see, e.g. \cite{Bomze2012, deKlerk,Dur2010}.

In {\it linear} CoP, the objective function is linear and the constraints are formulated with the help of linear matrix functions. Linear copositive problems are closely related to that of linear {\it Semi-Infinite Programming} (SIP) and
{\it Semidefinite Programming} (SDP). Copositive and semidefinite problems are particular cases of SIP problems, but CoP deals with  more challenging and less studied  problems than  SDP. The literature on the theory and methods of  SIP, CoP, and  SDP is quite extensive. We refer the interested readers to \cite{AH2013,HandbookSDPnew,Bomze2012, Dur2010,Weber2,W-handbook}, and the references in these works.

In   convex and conic  optimization, optimality conditions, and  duality results are usually formulated under  certain  regularity conditions, so-called constraint qualifications (CQ) (see, e.g. \cite{HandbookSDPnew,Kort,solodov,W-handbook}). Such conditions should guarantee the fulfillment of  the Karush-Kuhn-Tucker (KKT)- type optimality conditions and the {\it strong duality} property  consisting in the fact that the optimal values of the primal problem and  the corresponding  Lagrangian dual  one  are equal and the dual problem attains its maximum. Strong duality is the cornerstone of convex optimization, playing a particularly important role in the stability of numerical methods.

Unfortunately, even in  convex optimization, many  problems cannot be  classified  as regular  (i.e. satisfying some regularity conditions  such as, for example, strict feasibility). In  \cite{Dry-Wol}, we  read: ``...{\it new optimization modeling techniques and convex relaxations for hard nonconvex problems have shown that the loss of strict feasibility is a more pronounced phenomenon than has previously been realized}''.
 This phenomenon can  occur because of either   the poor choice of functions that describe feasible sets or the degeneration of the feasible sets themselves. According to \cite{Tuncel}, sometimes the loss of a certain CQ ``...{\it  is a modeling issue rather than inherent to the problem instance...}'' which ``... {\it justifies the pleasing paradigm: efficient modeling provides for a stable program}''.

Thus, the idea of a {\it regularization}  appears quite naturally which is aimed at obtaining an equivalent and more convenient  reformulation of the problem with some required  properties, one of which is that the regularized  problem  must satisfy  the generalized Slater condition.

The first papers on regularization of abstract convex problems (the regularization procedures are called {\it pre-processing} there)  appeared in the 1980-th, being followed by various  publications  for special classes of conic problems   (see, e.g. \cite{BW4,BW1}).  Nevertheless, as  Drusvyatskiy and  Wolkowicz wrote in \cite{Dry-Wol} published in 2017, for conic optimization in general, research in the field of regularization   algorithms is still in its infancy. At the same time, the authors of \cite{Dry-Wol} confirm that  to make a regularization algorithm viable, it is necessary to actively explore the structure of the problem since for some specific applications of conic optimization, the rich basic structure makes regularization quite possible and leads to significantly simplified models and enhanced algorithms.

Several approaches to the regularization of conic optimization problems are proposed in the literature.
In \cite{BW4, BW1}, the concept of the {\it minimal cone} of constraints was used by  Borwein and  Wolkowicz for regularization of abstract convex and conic convex  problems for which any CQ fails. An algorithm, proposed there for the description of the minimal cone, is based on a successive reduction of cone's faces  and  was named by the authors the {\it Facial Reduction Algorithm (FRA)}.

A different approach was proposed by  Luo,  Sturm, and  Zhang (see \cite{Luo} and the references therein) which is called  the {\it dual regularization}  or {\it conic expansion}. This approach tries to close a {\it duality gap} (the difference between the primal and dual optimal values) of the regularized problems  by expanding the  dual constraints' cone.

In \cite{Waki},  Waki and Muramatsu  applied the  facial reduction approach  to a conic optimization problem in such a way that each primal reduced cone is  dual to the cone generated by the conic expansion approach.

The facial reduction  approach has been successfully applied to  SDP and  the second-order cone programming problems, as well as to  certain classes of optimization problems over symmetric  (i.e. self-dual and homogeneous) and nice  cones (see, e.g. \cite{PP, Perm,Terlaki,Ramana-W}). At the same time, the question of effective constructive application of this approach to other classes of problems remains open. This is because   the known FRAs are more {\it conceptual }  than  practical.

In this paper, based on the results from  \cite{KT-new,KT-JOTA, KT-SetValued}, we develop a different approach to  regularization  of linear  CoP  problems. This approach is based on the concept of {\it immobile indices}, i.e. indices of the constraints that are active for all feasible solutions.

\vspace{2mm}

The purpose of the paper is to

\vspace{-2mm}

\begin{enumerate}

\item[a)]  describe in details a finite  algorithm for regularization of linear CoP problems that   is based on the concept of immobile indices but  does not require any additional information about  them;

  \vspace{-4mm}

  \item[ b) ] compare two approaches to regularization of  the  linear   CoP problems, one based on  the  facial reduction and the other based on the concept of immobile indices,   and the  corresponding regularized problems constructed using these approaches.

\end{enumerate}

\vspace{-2mm}

To the best of our knowledge, in CoP there has never  been  an attempt to develop detailed and easy in implementation algorithms, based on the minimal cone representation (see, e.g. the FRA  in  \cite{BW4,BW1} and the modified  FRA in \cite{Waki}). Nor do we have any information about   any other attempts to describe constructive regularization procedures for linear copositive problems. The regularization algorithms presented in the paper are new, original, and timely due to the growing number of eminent applications of CoP.

\vspace{2mm}

The paper is organized as follows. Section \ref{Introduction} hosts Introduction. Section \ref{optimality} contains  equivalent formulations of the linear CoP problem and the basic definitions. In section \ref{sec5}, we describe different   ways of  regularization of  copositive problems. In subsection \ref{subsec1}, we describe the minimal face regularization from \cite{ BW4, BW1} when it is  applied to  linear CoP problems; in  \ref{subsec2}, we shortly describe the proposed in \cite{KT-SetValued} one-step  regularization based on the concept of immobile indices and  compare the regularized problems obtained in subsections \ref{subsec1} and \ref{subsec2}. Section \ref{Sec6} contains iterative algorithms  for regularization of linear copositive problems: the   Waki and Muramatsu's facial reduction algorithm is described in \ref{subsec5}, a new regularization  algorithm based on the immobile index set  together with its compressed  modification  is  described in  \ref{subsec6}, followed by a short discussion on the iterative algorithms.  Section \ref{Concl} contains some conclusions.

\section{Linear copositive  programming  problem: equivalent \\ formulations and basic definitions}\label{optimality}

Given an integer $p>1$, denote by  $\mathbb R^{ p}_+$  the set of all $ p$ vectors with non-negative components,
by ${\mathcal S}(p)$ and $\mathcal S_+(p)$  the space of real symmetric $p\times p$ matrices
and the cone of symmetric positive semidefinite $p\times p$ matrices, respectively, and let
$\mathcal{COP}^{p}$  stay for the cone of symmetric copositive $p\times p$ matrices:
$$\mathcal{COP}^p:=\{D\in {\mathcal S}(p):t^{\top}Dt\geq 0 \ \forall t \in \mathbb
R^p_+\}.$$ The space $\mathcal S(p)$ is considered here as a vector space with the trace inner product $$A\bullet B:={\rm trace}\, (AB).$$

Consider a linear copositive programming  problem in the form

\begin{equation} \displaystyle\min_{x \in \mathbb R^n } \ c^{\top}x, \;\;\;
\mbox{s.t. } {\mathcal A}(x)
\in \mathcal{COP}^p,\label{SDP-2}\end{equation}
 where $x=(x_1,...,x_n)^{\top}$  is the vector of  decision variables. The data of the problem are presented by  vector $c \in \mathbb R^n$  and the constraints matrix function $\mathcal A(x)$  defined  in the form
\begin{equation}\label{A} \mathcal A(x):=\displaystyle\sum_{i=1}^{n}A_ix_i+A_0,\end{equation}
 with given  matrices $A_i\in \mathcal S(p), i=0,1,\dots,n$.

It is well known (see e.g. \cite{AH2013}) that the  copositive problem   (\ref{SDP-2}) is equivalent to the following convex  SIP problem:
\begin{equation} \displaystyle\min_{x \in \mathbb R^n } c^{\top}x, \;\;\;
\mbox{s.t. } t^{\top}{\mathcal A}(x)t
\geq  0 \;\; \forall t \in T,\label{SIP}\end{equation}
 with a $p$ - dimensional compact index set in the form of a simplex
 \begin{equation}\label{SetT}T:=\{t \in \mathbb R^p_+:\mathbf{e}^{\top}t=1\},\end{equation}
where $\mathbf{e}=(1,1,...,1)^{\top}\in \mathbb R^p. $

Denote by $X$ the  feasible set  of the equivalent problems (\ref{SDP-2}) and   (\ref{SIP}):
\be X:=\{x\in \mathbb R^n: {\mathcal A}(x)\in {\cal COP}^p\}
=\{x\in \mathbb R^n:t^{\top}{\mathcal A}(x)t\geq 0 \;\; \forall  t \in T\}.\label{setX}\ee

In what follows, we will suppose that $X\not=\emptyset.$

Evidently,  the set $X$ is convex.

\begin{remark}\label{R-1}
 Since $X\not=\emptyset,$ then without loss of generality, we can consider that $A_0\in {\cal COP}^p. $ Indeed, by fixing a feasible solution $y\in X$ and substituting  the variable  $x$  by  a new  variable $z{:}=x-y$, we can replace the original problem
(\ref{SDP-2}) by the following one in terms of $z$:

\vspace{-3mm}

\begin{equation*} \displaystyle\min_{z \in \mathbb R^n} \ c^{\top}z, \;\;\;
\mbox{s.t. } \bar {\mathcal A}(z)
\in \mathcal{COP}^p,\end{equation*}

\vspace{-3mm}

with
$\bar {\mathcal A}(z):=\displaystyle\sum_{i=1}^{n}A_iz_i+\bar A_0,$ $ \bar A_0:= \mathcal A(y)\in {\cal COP}^p.$
\end{remark}

 According to the commonly used  definition, the constraints of the copositive problem (\ref{SDP-2}) satisfy {\it the  Slater condition} if
\begin{equation}\label{Slait1}
\exists \; \bar{x} \in \mathbb R^n \; \mbox{ such that } \; {\mathcal A}(\bar x) \in {\rm int }\, \mathcal{COP}^p=
\{D\in \mathcal{S}(p): t^{\top}Dt>0 \; \forall t \in \mathbb R^p_+,\; t \not=\bf{0}\}.
\end{equation}

Here ${\rm int  }\,\mathcal B$   stays for the interior  of a set $\mathcal B$.

\vspace{2mm}

Following \cite{KT-JOTA, KT-SetValued}, let's define the {\it set of normalized  immobile indices} $ T_{im}$  in  problem (\ref{SDP-2}):

\vspace{-3mm}

\begin{equation}\label{Tnormalized}  T_{im}:=\{t \in   T : t^{\top} \mathcal A (x)t=0 \ \  \forall x\in {X}\} .\end{equation}
In what follows, the elements of the set $T_{im}$ are called immobile  indices.


The following lemma follows from  Lemma 1 and  Proposition 1 in \cite{KT-SetValued}.
\begin{lemma} \label{LOK} Given the linear copositive problem (\ref{SDP-2}),   (i) the Slater condition (\ref{Slait1}) is equivalent to the emptiness of the  set $T_{im},$
(ii) the normalized immobile index set $
T_{im}$ is either empty or  can  be represented as a union of a finite number of convex  closed  bounded polyhedra.
\end{lemma}

 For a vector  $t=(t_k,k\in P)^{\top}\in \mathbb R^p_+$ with   $ P:=\{1,2,...,p\}$, define the sets
$$P_+(t):=\{k\in P:t_k>0\}, \ P_0(t):=P\setminus P_+(t).$$
 Given a set $\mathcal B$  and a point $l=(l_k,k \in P)^{\top}$  in $\mathbb R^p$,  denote by $\rho(l,\mathcal B)$ the distance between these set and point,
$\ \rho(l,\mathcal B):=\min\limits_{\tau \in \mathcal B}\sum\limits_{k\in P}|l_k-\tau_k|, $ and  by
${\rm conv } \mathcal B$   the convex hull of the  set $\mathcal B$.

\vspace{3mm}

 Suppose that in problem (\ref{SDP-2}), the normalized  immobile index set $T_{im}$ is non-empty.
Consider a finite  non-empty subset of   $T_{im}$:
\be V=\{\tau(i)\in T_{im} \; \forall i \in I\},\; 0<|I|<\infty.\label{N0**}\ee
    For  this set, define the following number and sets:
\be \sigma(V):=\min\{\tau_k(i),  \ \  k \in P_+(\tau(i)) , \; i \in I\}>0,\ \ \qquad\qquad\qquad\qquad\label{ep}\ee
\be \label{omega}\Omega(V):=\{t \in T:\rho(t,{\rm conv} V)\geq \sigma(V)\},\qquad\qquad\qquad\qquad\qquad\qquad\ee
\begin{equation} \mathcal X(V):=\{x \in \mathbb R^n:
{\mathcal A}(x) \tau(i)\geq 0 \; \forall i \in {I}; \;\;  t^\top {\mathcal A}(x)t\geq 0 \;  \forall t \in \Omega(V)\}. \label{calX}\end{equation}

In  \cite{KT-new}, the following theorem is proved.

\begin{theorem}\label{L25-02-1}  Consider problem (\ref{SDP-2}) with  the  feasible set $X$.  For any subset  (\ref{N0**})   of the set of normalized immobile indices of this problem, the following {equality} holds true:
$$ X={\mathcal X}(V),$$
where the set  ${\mathcal X}(V)$ is defined in (\ref{calX}).
\end{theorem}

   \section{Regularization of copositive problems}\label{sec5}

In this section,  first, we  remind a known regularization approach  developed in \cite{ BW4, BW1} for conic optimization problems and based on the concept of the minimal face.  We  briefly  describe how   this approach can be applied to linear CoP problems. After, for  the copositive problem (\ref{SDP-2}), we  present  another  regularization approach  based on the concept of immobile indices  and compare the regularized problems obtained using two considered approaches.

\vspace{2mm}

  \subsection{Minimal face regularization}\label{subsec1}

Let us, first, recall  the necessary terms and notions.

	\vspace{2mm}
	
By definition, a  convex subset $\mathbf{F}$
   of the cone ${\cal COP}^p$ is its {\it face}   if for any $x\in {\cal COP}^p$, $y\in {\cal COP}^p$,
  the inclusion  $x + y\in \mathbf{F}$ implies $x\in \mathbf{F},$ $y\in \mathbf{F}.$ It is evident that any face of the cone ${\cal COP}^p$ is also a cone.

Given the copositive problem (\ref{SDP-2}) with the feasible set $X$   presented  in (\ref{setX}), let $\mathbf{F}_{min}$ be the
smallest (by inclusion) face of ${\cal COP}^p$  containing a  set  $\cal D $  defined  in terms of the constraints of this problem  as follows:
\be {\cal D}:=\{ \mathcal A(x), \ x \in X\}.\label{D}\ee

 In what follows, the face  $\mathbf{F}_{min}$ will be called {\it the minimal face of the optimization problem } (\ref{SDP-2}).

	\vspace{3mm}
	
 Generally speaking,  for the copositive  problem (\ref{SDP-2}), the approach suggested in \cite{ BW4, BW1},   is to replace the constraint ${\cal A}(x)\in {\cal COP}^p$ with an equivalent constraint ${\cal A}(x)\in \mathbf  F_{min}$.   The resulting regularized problem    takes the form
\be \min_{ x\in \mathbb R^n }\ c^\top x, \;\; \mbox{ s.t. } {\mathcal A}(x)\in \mathbf{F}_{min}.\label{25P1}\ee

The dual problem to (\ref{25P1})  can be written in the form
\be \max_{ U \in \mathcal S(p) }\; -A_0\bullet U,\; \mbox{ s.t. } A_j\bullet U=c_j \; \forall j=1,...,n;\; U\in \mathbf{F}^*_{min},\label{25D1}\ee
where $\mathbf{F}^*_{min}$ is the dual cone to the cone $ \mathbf{F}_{min}.$

\vspace{3mm}

It is proved in \cite{ BW4, BW1}, that   the constraints of problem (\ref{25P1}) satisfy the {\it generalized Slater condition}:  there exists $\bar x\in X$ such that ${\mathcal A}(\bar x)\in {\rm relint} \; {\mathbf{F}}_{min}$ and hence  the duality gap between the dual problems (\ref{25P1}) and (\ref{25D1}) vanishes. Here ${\rm relint }\,\mathcal B$ stays for the relative interior of a set $\mathcal B$.

\vspace{3mm}

Unfortunately,  there is no available   information about how to explicitly construct the cones $ \mathbf{F}_{min}$ and $\mathbf{F}^*_{min}$ in a general case  and, in particular, in the case of copositive problems.

\subsection{One-step  regularization based on the concept of immobile indices}\label{subsec2}

In our paper  \cite{KT-SetValued},  for the  copositive problem (\ref{SDP-2}), we  obtained a  regularized dual  problem  which is different from (\ref{25D1}) . The construction of  this dual  is  based on
  the concept of immobile indices and  can be considered as  {\it one-step  regularization}  since it  contains  a unique   step.

\vspace{3mm}

Consider the copositive  problem  (\ref{SDP-2}). Let $T_{im}$ be the  defined in (\ref{Tnormalized})  normalized set of immobile indices  of this problem. If $T_{im}=\emptyset$, then  problem  (\ref{SDP-2}) satisfies the Slater condition,   which means that it  is already regular and no regularization is required.

\vspace{2mm}

Now, suppose that  $T_{im}\not =\emptyset$. In this case, the Slater condition is not satisfied and the problem  is not regular.  Let us describe  how one can convert    problem (\ref{SDP-2}) into a regularized one.

\vspace{2mm}

 Consider  the  set    ${\rm conv}\, T_{im}$ and the set  $W$ of  all  vertices  of ${\rm conv}\, T_{im}$:
\be W:=\{t(j), \; j \in J\}, \; 0<|J|<\infty.\label{ver}\ee
 Suppose that the elements $t(j), \; j \in J,$ of the set $W$ are known. Then we can regularize problem (\ref{SDP-2}) in just \textbf{one step}.

\vspace{3mm}

In fact, it follows from Theorem \ref{L25-02-1} and  Lemmas 2 and 3  in \cite{KT-SetValued} that the set $X$ of feasible solutions of the original problem (\ref{SDP-2})
coincides with the set of feasible solutions of the following system:
$$ t^\top{\mathcal A}(x)t\geq 0\; \forall  t \in \Omega(W);\;\;\;
\mathcal A(x)t(i)\geq 0 \; \forall i \in J,$$
and   the next condition is satisfied:
\be \exists \; \bar{x} \in X \; \mbox{ such that } \;  t^\top{\mathcal A}(\bar x)t>0\;\; \forall  t \in \Omega(W).\label{S-type}\ee
Here the set $\Omega(W)$  is   defined by the rules (\ref{omega}) with $V=W.$

Consequently, the original  copositive  problem (\ref{SDP-2}) is equivalent to the following  SIP problem:
 \bea & \min\limits_{x\in \mathbb R^n} \; c^\top x,\qquad \label{reg-1}\\
 &\mbox{s.t. }\;\; t^\top{\mathcal A}(x)t\geq 0\  \forall  t \in \Omega(W),\label{reg-2}\\
&{\mathcal A}(x)t(i)\geq 0\  \forall  i \in J.\label{reg-3}\eea
Problem (\ref{reg-1})-(\ref{reg-3}) can be considered as a {\it regularized} primal  problem since

\begin{itemize}

\vspace{-2mm}

\item  it possesses  a finite number of linear inequality constraints (\ref{reg-3}),

\vspace{-3mm}

\item the first group of  constraints, (\ref{reg-2}), satisfies the Slater type condition (\ref{S-type}),

\vspace{-2mm}

\item  the set $\Omega(W)$ is compact.
\end{itemize}

 Let us stress that in problem  (\ref{reg-1})-(\ref{reg-3}),   the infinite  index set   $\Omega(W)$  is obtained by removing the set $T_{im}$ together with the $ \sigma(W) $-neighborhood
 of its convex hull, from the original index set   $T$. Note here that  the set $\Omega (W)$
 \begin{itemize}
 \item[(a)] is  explicitly constructed  by the   rules  (\ref{ep}),  (\ref{omega}),  using the finite  set $W=\{t(j), j \in J\}$ of vertices of ${\rm conv} T_{im}$,

 \vspace{-3mm}

 \item[(b)] does not contain the set ${\rm conv}\,T_{im}$,

 \vspace{-3mm}

 \item[(c)]  may be sufficiently  small.
 \end{itemize}
 All these   properties may be useful for numerical  solving the  problem (\ref{reg-1})-(\ref{reg-3}).

\vspace{3mm}

It is evident that problem (\ref{reg-1})-(\ref{reg-3}) can be written in the equivalent conic form
\be \min_{ x\in \mathbb R^n }\ c^\top x,\;\; \mbox{ s.t. }
{\mathcal A}(x)\in {\cal K}_0,\label{25P2-2}\ee
where $ \  {\mathcal K}_0:=\{D\in {\cal S}(p):\ t^\top Dt\geq 0\;\;\forall t\in  \Omega(W); \  Dt(j)\geq 0 \; \forall j \in J \}.$

  It can be shown   that ${\cal K}_0\subset {\cal COP}^p$.

  The dual  problem  to (\ref{25P2-2})   is as follows:
\be \max_{ U \in \mathcal S(p)}\; -A_0\bullet U,\; \mbox{ s.t. } A_j\bullet U=c_j\; \forall j=1,...,n;\; U\in {\mathcal K}^*_{0}.\label{25D2}\ee
  In the problem above,  ${\mathcal K}^*_{0}$  is the dual cone to  ${\mathcal K}_0$   and has the form
\be  {\mathcal K}^*_{0}={\rm cl} \{D \ \in  {\cal S}(p): \ D{\in}{\cal CP}(W)\oplus {\cal P}^*\},\label{K**}\ee
where

\vspace{-5mm}

\bea& {\cal CP}(W):={\rm conv} \{tt^{\top}:t \in \Omega(W)\},\;\nonumber\\ & \  {\cal P}^*:=\{D \in  {\cal S}(p): \ D= \sum\limits_{j\in J}(\lambda(j)(t(j))^\top+t(j)(\lambda(j))^\top),
\ \lambda(j)\geq 0\; \forall  j \in J\}.\nonumber\eea

Here  and in  what follows, for  given sets $ {\mathcal B}$ and ${\mathcal G}$,  ${\rm cl} \, {\mathcal B}$ denotes the closure of the set ${\mathcal B}$ and ${\mathcal B}\oplus {\mathcal G}$ stays for the Minkowski sum of the  corresponding  two sets.

\vspace{2mm}

 Notice that  for the pair of dual  conic  problems (\ref{25P2-2}) and (\ref{25D2}), the duality gap is zero.

\vspace{2mm}

As  it was shown in \cite{KT-SetValued}, the cone (\ref{K**}) in problem (\ref{25D2}) can be replaced by the following one (which has a more explicit  form since it does not contain the   closure operator):
$${\mathbf{{K}}_0^*:}= \{D\in{\mathcal S}(p): \ D\in {\cal CP}^{p}\oplus  {\cal P}^*\},$$
where $\mathcal{ CP}^{p}$  denotes the set of completely positive matrices:
\be {{\cal CP}^{p}}:={\rm conv} \{tt^{\top}:t \in \mathbb R_+^{p}\},\label{CP}\ee
and there is no  duality gap for  problem (\ref{25P2-2}) and  its dual problem in the form (\ref{25D2})  with ${\mathcal K}^*_{0}$ replaced by  $ {\mathbf{{K}}_0^*}$.

\vspace{3mm}

Note that  the cones ${\mathcal K}_0$ and ${\mathbf{{K}}_0^*}$ are  explicitly described in terms of  indices (\ref{ver}) and this is an advantage of the presented here approach over  that  described in \ref{subsec1}.

The only drawback of  the   regularization  procedure  described here is  the following:
\begin{itemize}
\item []  {\it to apply   the one-step regularization , one needs to know the  finite number of indices (\ref{ver})
  which are the vertices of the set ${\rm conv}\, T_{im}.$}
\end{itemize}

 It is easy to see  that the  regularized primal problem (\ref{25P2-2}) can be modified as follows:
\be  \min\limits_{x\in \mathbb R^n} \; c^\top x,\;\; \mbox{ s.t. } {\mathcal A}(x)\in \overline {\mathcal K}_0,\label{R-2}\ee
where
\begin{equation*}\begin{split} \overline {\mathcal K}_0=\{D\in S(p): \ & t^\top Dt\geq 0\;\; \forall t \in \Omega(W), \\
&e^\top_kDt(j)=0 \; \forall k \in P_+( t(j) ),\;\;  e^\top_k Dt(j)\geq 0 \; \forall k \in P\setminus P_+( t(j) ),\ \forall j \in J\},\label{barK*}\end{split}\end{equation*}
 $\{e_k, \; k \in P\}$ is the standard (canonical) basis of the vector space $\mathbb R^p$.

\vspace{3mm}

It is evident that $\overline {\mathcal K}_0\subset  {\mathcal K}_0 $ and, as  mentioned above,   ${\cal K}_0\subset {\cal COP}^{p}$. Hence $\overline {\cal K}_0\subset {\cal COP}^p.$

 To show that  these regularizations are themselves deeply connected,  let us give  an  explicit description of the minimal face $ \mathbf{F}_{min}$ in terms of the vertices  of the set ${\rm conv} \, T_{im}$ and the index sets $M(j),\ j \in J,$  defined  as:
 \be M(j): =\{k \in P:e^\top_k{\cal A}(x)t(j)=0\;\; \forall x \in X\}, \   j \in J.\label{N*}\ee

 The following theorem can be proved (see  \cite{KT-new}).

\begin{theorem}\label{T3}    Given copositive problem (\ref{SDP-2}), let $\{t(j),j\in J\}$ be the (finite) set of all  vertices  of the set ${\rm conv}\, T_{im}$. Then the minimal face $ {\mathbf{F}}_{min}$ of this problem can be described in two equivalent forms
\begin{equation*}\begin{split}\mathbf{F}_{min}=K_{min}(1):=\{D\in {\cal COP}^p: & \ e^\top_kDt(j)=0 \  {\forall} k \in {M(j)}, \;  \forall  j \in J\}, \ \mbox{and}\\
\mathbf{F}_{min}={K}_{min}(2){:}=\{D\in {\cal COP}^p: &\ e^\top_kDt(j)=0 \ {\forall} k \in {M(j)},\\
&\ e^\top_kDt(j)\geq 0 \  \forall  k \in {P\setminus M(j)},  \forall  j \in J\}.\end{split}\end{equation*}
\end{theorem}

Now, having described the minimal  face  $\mathbf{F}_{min}$ via immobile indices, we can compare  the  regularized problems (\ref{25P1}), (\ref{25P2-2}), and  (\ref{R-2})  in more detail.

 The regularized problem (\ref{25P1}) is formulated  using  the facial reduction approach  to the copositive problem (\ref{SDP-2}) and  the  regularized problems (\ref{25P2-2}) and  (\ref{R-2}) are  obtained using  the immobile indices of this problem. The difference between these three problems  is    that in  problem (\ref{25P1}), the constraint set is determined by  the minimal face  $\mathbf{F}_{min}$, while the constraints of  problem  (\ref{25P2-2}) are formulated with the help of  the cone $\mathcal K_0$, and   the constraints of problem (\ref{R-2}) use the cone $\overline {\mathcal K}_0$.
It should be noticed that the minimal face  $\mathbf{F}_{min}$ and the  cones  $\mathcal K_0$ and $\overline {\mathcal K}_0$ satisfy  the inclusions $$\mathbf{F}_{min}\subset \overline {\mathcal K}_{0}\subset {\mathcal K}_{0}.$$ At the same time, the cones $\mathbf{F}_{min}$ and
{$\overline {\mathcal K}_{0}$} are faces of the  cone of copositive matrices ${\cal COP}^{p}$, while the  cone ${\mathcal K}_{0}$ is not,  in general. One can show that {$\overline {\mathcal K}_{0}$} is an exposed face while the face $\mathbf{F}_{min}$ is not, in general.

\vspace{3mm}

  For each of the mentioned above conic  problems, we face certain challenges caused by the  troubles  connected with the  {\it concrete construction} of the respective cones.
 Thus, for example, for the copositive problem (\ref{SDP-2}),  the following difficulties should be mentioned:
\begin{itemize}

\item to define the cones ${\mathcal K}_0$ and  $\overline {\mathcal K}_0$, the elements  $t(j), \; j \in J,$  of the finite set of indices (\ref{ver})  should be known;

\vspace{-3mm}

\item as far as we know, there are no explicit procedures of  constructing the minimal face  $\mathbf{F}_{min}$
 and its dual cone $\mathbf{F}^*_{min}$.

\end{itemize}

 Theorem \ref{T3} shows  how the minimal face $\mathbf{F}_{min}$ can be represented in the  form of the cones ${K}_{min}(1)$ and $K_{min}(2)$ via immobile indices.  Notice that to construct these cones, one has to find not only the set of indices (\ref{ver}) but also the corresponding  sets $ M (j), j \in J$, defined in (\ref{N*}).

\vspace{3mm}

As mentioned above,  regularity is an important property of optimization problems. As a rule, the  regularity of copositive problems is characterized by the Slater condition. In this regard, it is important to note that the  regularized problem (\ref{25P1}) satisfies the {\it generalized } Slater condition while the obtained here regularized problems (\ref{25P2-2}) and (\ref{R-2}) satisfy the  Slater type condition (\ref{S-type}).  This difference can be  important for further study of linear CoP problems, as well as for the development of stable numerical methods for them.

\section{Iterative algorithms  for regularization of linear copositive problems}\label{Sec6}

In section \ref{sec5}, we considered general schemes of two  theoretical methods  that allowed us to obtain  regularizations of the  linear  copositive problem (\ref{SDP-2}). In each of these {schemes},  we  meet  some difficulties associated with explicit representations of the  respective "regularized" feasible cones and their dual ones.
In this section, we consider and compare two different regularization approaches aimed to overcome these difficulties using algorithmic procedures.

\subsection{  Waki and Muramatsu's facial reduction algorithm}\label{subsec5}

 In \cite{Waki},  for linear conic problems,  a  regularization   algorithm  was proposed  by  Waki and  Muramatsu. This algorithm can be considered as the Facial Reduction  Algorithm (FRA) from  \cite{BW4, BW1} applied to linear  conic problems in finite-dimensional spaces.

\vspace{2mm}

Let us    describe the algorithm from \cite{Waki} for the linear copositive  problem (\ref{SDP-2})  with the matrix constraint function $\mathcal A(x)$ defined in (\ref{A}). Remind that  here  we  consider that  problem (\ref{SDP-2})  is feasible. Then, according to Remark \ref{R-1},  we  can  assume  that $A_0\in {\cal COP}^{p}$.

Denote
$${\rm Ker}\,\mathcal A:=\{D\in S(p): A_j\bullet D=0 \ \forall j=0,1,...,n\}.$$
 As above,  let ${\mathcal F}^*$ denote the dual cone of a given cone ${\mathcal F}$.

\vspace{3mm}

  For a given  feasible copositive problem (\ref{SDP-2}), starting  with ${\cal COP}^{p}$, the  Waki and Muramatsu's algorithm repeatedly finds smaller faces of ${\cal COP}^{p}$ until it  stops with   the minimal face $\mathbf{F}_{min}$.

\vspace{3mm}

{\bf   Waki and Muramatsu's  FRA   for the copositive problem (\ref{SDP-2})}

\vspace{3mm}

$\qquad$ {\bf Step 1:} Set $i{:}=0$ and ${\mathcal F}_0:={\cal COP}^{p}.$

$\qquad$ {\bf Step 2:} If ${\rm Ker}\mathcal A\cap {\mathcal F}^*_i\subset {\rm span}\{Y_1,...,Y_i\},$ then STOP:
 $\mathbf{F}_{min}={\mathcal F}_i.$

$\qquad$ {\bf Step 3:} Find $Y_{i+1}\in {\rm Ker}\mathcal A\cap {\mathcal F}^*_i\setminus  {\rm span}\{Y_1,...,Y_i\}.$

$\qquad$ {\bf Step 4:} Set ${\mathcal F}_{i+1}{:}={\mathcal F}_i\cap \{Y_{i+1}\}^\bot $  and $i:=i+1$,  and go to step 2.

\vspace{3mm}

 The description of the algorithm is very simple  but, in practice, its  implementation presents serious difficulties which arise on   step 2 and  especially on  step 3.  As the matter of fact, in the case of the copositive problem (\ref{SDP-2}), the fulfillment  of  step 3 is hard  already on the first two iterations.

  Let us consider the initial iteration when $i=0$.  On  step 3, one has to find a matrix
 $Y_1\in {\rm Ker}\mathcal A\cap {\mathcal F}^*_0.$
Since ${\mathcal F}_0={\cal COP}^{p}$, then  at  the current  iteration  ($i=0$) we know the explicit description of the dual cone for ${\mathcal F}_0$:
${\mathcal F}^*_0={{\cal CP}^{p}}$, where  the cone  ${\cal CP}^{p}$ is  defined in (\ref{CP}). Therefore, the matrix $Y_1$ should have the form
\begin{equation*} Y_1=\sum\limits_{i\in I_1}t(i)(t(i))^\top, \; t(i)\geq 0,\;  t(i)\not= 0 \ \forall i \in I_1,\;\; 0<|I_1|\leq p(p+1)/2,\label{3q}\end{equation*}
and  the condition
$ \sum\limits_{i\in I_1}(t(i))^\top A_jt(i)=0\;\forall j=0,1,...,n$,  has to  be satisfied. 

 At  the next iteration ($i=1$), one is looking for a matrix $Y_2$ such that
\begin{enumerate}
\item[\bf{C1}:]  $Y_2\in {\mathcal F}^*_1={\rm cl} \{D\in S(p):D\in{\cal CP}^{p}\,\oplus\; \alpha Y_1,\; \alpha  \in \mathbb R\},$

    \vspace{-2mm}

\item[\bf{C2}{:}]  $Y_2\not \in {\rm span}\{Y_1\},$

\vspace{-2mm}

\item[\bf{C3}{:}]  $A_j\bullet Y_2=0 \  \forall   j=0,1,...,n.$
\end{enumerate}

The first difficulty arises when trying to satisfy the  condition \textbf{C1}, as there is no explicit description of the set
${\mathcal F}^*_1$.  Notice that this set is defined using the closure operator, this operator being  essential for the definition of ${\mathcal F}^*_1$.
Therefore, in general, for  a matrix  $Y_2$ satisfying the condition {\bf C1}, it may happen that
$Y_2\not \in \{D\in S(p): \ D\in {\cal CP}^{p}\,\oplus\; \alpha Y_1, \; \alpha\in \mathbb R\}.$

\vspace{2mm}

In \cite{Waki},  there is {also} no any indication of how to find a matrix $Y_2$ satisfying the conditions { \textbf{C2}} and { \textbf{C3}}.  Notice that {the} fulfillment of   these  conditions  is a non-trivial task  as well.

\vspace{3mm}

Thus,  we can state that although  the  reported in \cite{Waki}  FRA   is  an easy-to-describe method,  its practical implementation is not constructively  described, which makes it  difficult to apply.
There is no information concerning which form should have the matrix  $Y_i$ at the $i$-th iteration ( $i\geq 1$) of the algorithm and how  to meet   the conditions  \textbf{ C1} - { \textbf{C3}}  for it.

\subsection{A regularization   based on the immobile indices }\label{subsec6}

Here we will describe and justify a different  algorithm  for  regularization of the copositive problem (\ref{SDP-2}). This algorithm has a similar structure to the  Waki and  Muramatsu's FRA   considered   in \ref{subsec5} but  is based on the concept of immobile indices and  described in more detail,  being, therefore,  more constructive. Note from the outset that although   our  algorithm exploits the properties of the set of immobile indices, it does not require the initial knowledge of either this set or the vertices of its convex hull.

\subsubsection{Algorithm \textbf{REG-LCoP} (REGularization of Linear Copositive Problems)}\label{4-2-1}

\vspace{2mm}

{\bf Iteration $\#$ 0.} Given the copositive  problem in the form (\ref{SDP-2}), consider the following {\it regular} SIP  problem:
\begin{equation*}{SIP}_0: \qquad \displaystyle\min_{(x,\mu)\in \mathbb R^{n+1}}\ \mu, \mbox{ s.t. }
t^{\top}{\mathcal A}(x)t+\mu\geq 0\;  \forall  t \in T,\label{1a}\end{equation*}
with the index set $T$ defined in (\ref{SetT}).

If there exists a feasible solution $(\bar x,\bar \mu)$ of this problem with $\bar \mu<0$, then set
$m_*:=0$   and  go to the {\it Final step.}

 Otherwise  the vector $(x=\mathbf 0,\mu=0)$  is an optimal solution of the problem ($SIP_0$).

It should be noticed that in the problem ($SIP_0$), the index set $T$ is compact, and the constraints  satisfy the
Slater condition.  Hence, (see e.g. \cite{Bon}),  it follows from the optimality conditions for  the vector $(x= \mathbf 0 ,\mu=0)$ that there exist indices and numbers
    \begin{equation*}\tau(i)\in T, \; \; \gamma(i)>0  \ \forall i \in I_1, \;\; |I_1|\leq n+1,\end{equation*}

such that
$\sum\limits_{i\in I_1}\gamma(i)(\tau(i))^{\top}
A_j\tau(i)=0  \  \forall   j=0,1,...,n;\;\;  \sum\limits_{i\in I_1}\gamma(i)=1.$

It follows from {the relations above} that $I_1\neq \emptyset $ and $\tau(i)\in T_{im}\subset T \;  \forall  i \in I_1,$
  and hence $ e^\top_k{\cal A}(x)\tau(i)=0\;  \forall  k \in P_+(t(i)),\; \forall   i \in I_1,\  \forall x\in X.$ Set   $L_1(i):=P_+(\tau(i)), \; i \in I_1,$ and go to the next iteration.

\vspace{2mm}

{\bf Iteration $\#$ $m$, $m\geq 1$.}
By the beginning of  the iteration, we have indices and sets $\tau(i),$  $ L_m(i),$ $i\in I_m,$ such that
\be \tau(i)\in T_{im}, \ P_+(\tau(i))\subset  {L_m(i)}\subset P,\;  e^\top_k{\cal A}(x)\tau(i)=0\ \forall  k \in L_m(i),\;  \forall x\in X,\; \forall  i \in I_m.\label{vvv}\ee

Consider a SIP problem
\begin{equation*}\begin{split} & \min\limits_{(x,\mu)\in \mathbb R^{n+{1}}} \mu, \\
{SIP}_m: \qquad\mbox{ s.t. } \ \ &e^{{\top}}_k{\mathcal A}(x)\tau(i)
=0 \;  \forall  k \in L_m(i) ;\
e^{{\top}}_k{\mathcal A}(x)\tau(i)\geq 0 \;  \forall  k\in P\setminus L_m(i),\
 \; \forall  i \in I_m,\\
&t^{\top}{\mathcal A}(x)t+\mu\geq 0 \;  \forall  t \in  \Omega(W_m),\end{split}\end{equation*}
where $W_m:=\{\tau(i), i \in I_m\}$,  the set  $\Omega(W_m)$ is constructed by the rules  (\ref{ep}),  (\ref{omega}) with $V=W_m$.
Since  $W_m$ is  a subset   of the set of immobile indices
 of problem (\ref{SDP-2}), then it follows from Theorem
 \ref{L25-02-1}  and  the  equalities in (\ref{vvv}) that   ${\cal X}(W_m)=X$.

 Notice that by the definition of the set $\Omega(W_m)$, it holds:
$$\rho(t,{\rm conv} W_m)\geq \sigma(W_m) >0\;\; \forall t \in  \Omega(W_m).$$

\vspace{2mm}

In the  problem (${SIP}_m$),
the index set $\Omega(W_m)$ is   compact and the
constraints  satisfy
the following Slater type condition:
\begin{equation*}\begin{split} \exists (\widehat x,\widehat \mu) \mbox{ such that } & e^{{\top}}_k{\mathcal A}(\widehat x)\tau(i)
=0 \;  \forall  k \in {L_m(i)};\,
e^{{\top}}_k{\mathcal A}(\widehat x)\tau(i)\geq 0 \;  \forall k \in P\setminus {L_m(i)}, \;  \forall  i \in I_m, \\
 &t^{\top}{\mathcal A}(\widehat x)t+\widehat \mu> 0 \;  \forall  t \in \Omega(W_m).\end{split}\end{equation*}
Hence,  this problem  is {\it regular}.

\vspace{2mm}

If  problem {(${SIP}_m$)} admits a feasible solution $(\bar x,\bar \mu)$ with $\bar \mu<0$, then STOP and go to the
  {\it Final step}  with $m_*:=m.$

Otherwise, the   vector $(x=\mathbf 0,\mu=0)$  is an optimal solution  of {(${SIP}_m$)}.  Since this  problem is regular, the optimality of the vector $(x=\mathbf 0,\mu=0)$ provides (see \cite{levin}) that
there exist   indices,  numbers, and vectors
\be  \tau(i)\in \Omega(W_m), \;   \gamma(i), i \in \Delta I_m, \;   \; 1\leq |\Delta I_m|\leq n
+1;\   \lambda^m(i)\in \mathbb R^p,  i\in I_m,\label{11.1}\ee
which satisfy the following conditions:
\begin{equation}\begin{split} &\sum\limits_{i\in \Delta I_m}\gamma(i)(\tau(i))^{\top}A_j\tau(i)+
2\sum\limits_{i\in I_m}( \lambda^m(i))^{\top}A_j\tau(i)=0 \;  \forall
j=0,1,...,n;\\
&\gamma(i)>0 \;  \forall i \in \Delta I_m;\;\;\;  \lambda^m_k(i)\geq 0 \;  \forall  k \in P\setminus  {L_m(i)},\;  \forall  i \in I_m.\label{18-11}\end{split}\end{equation}
Here and in what follows, without loss of generality, we suppose that $\Delta I_m\cap I_m=\emptyset.$
Moreover, applying  the  described in {\cite{KT-Global}}  procedure \textbf{DAM} to the data (\ref{11.1}), it is possible to ensure that the following conditions are met:
\be P_0(\tau(i))\cap P_+(\tau(j))\not=\emptyset \; \; \forall i \in \Delta I_m, \; \forall j\in I_m.\label{11*}\ee

Let us set $\Delta L(i):=\{k \in P\setminus L_m(i):\lambda^m_k (i) >0\}, \  i \in I_m,$
$${L_{m+1}(i)}:={L_m(i)}\cup \Delta {L}(i), \  i \in I_m;\;\; {L_{m+1}(i)}:=P_+(\tau(i)), \  i \in \Delta I_m.\qquad$$
    It follows from (\ref{18-11}) that $e^\top_k{\cal A}(x)\tau(i)=0\  \forall   k \in \Delta L(i),$ $\ \forall   i \in I_m,$ and
 $\tau(i)\in T_{im}\; \forall   i \in \Delta I_m.$ The last inclusions imply the equalities $e^\top_k{\cal A}(x)\tau(i)=0\; \forall   k \in P_+(\tau(i)),$ $ \forall  i \in \Delta I_m.$

 Go to the next iteration $\# (m+1)$ with the new data
\be \tau(i),\; L_{m+1}(i), \; i \in I_{m+1}:=I_m\cup \Delta I_m.\label{o-18}\ee

\vspace{1mm}

{\bf Final step.}
It follows from (\ref{11*}) that the algorithm \textbf{\textbf{REG-LCoP}}  runs a finite number $m_*$ of iterations.  Therefore, for some $m_*\geq 0$,  the problem  (${SIP}_{m_*})$  has a feasible solution $(\bar x, \bar \mu)$ with $\bar \mu<0$.
 Observe  that   by  Theorem \ref{L25-02-1},    the found vector $\bar x$ is a feasible solution of the original copositive problem (\ref{SDP-2}).

 If $m_*=0$, then  the constraints of  problem (\ref{SDP-2}) satisfy the Slater condition with $\bar x$, and hence the problem is regular.

Suppose now that  $m_*>0$.
Consider   problem
\begin{equation*}\begin{split}& \min\limits_{x\in \mathbb R^{n}} c^\top x, \\
\mbox{ s.t. }& \ e^{{\top}}_k{\mathcal A}(x)\tau(i)=0 \;  \forall k \in {L_{m_*}}(i) ,\
e^{{\top}}_k{\mathcal A}(x)\tau(i)\geq 0\ \;  \forall  k\in P\setminus {L_{m_*}}(i),\;  \forall  i \in I_{m_*}, \\
&t^{\top}{\mathcal A}(x)t\geq 0 \;  \forall  t \in  \Omega(W_{m_*}),\end{split}\end{equation*}
where  the sets  $W_{m_*}=\{\tau(i), i \in I_{m_*}\}$  and $\Omega(W_{m_*})$    are  the same as in the problem (${SIP}_{m_*})$. Note that by construction,  ${\cal X}(W_{m_*})=X. $
  Hence  the problem above  is equivalent to  problem (\ref{SDP-2}) and  can be considered as
its regularization since

\vspace{-2mm}

\begin{itemize}

\item   it has a finite number of linear equality/inequality constraints,

\vspace{-2mm}

\item    by construction,  $t^{\top}{\mathcal A}(\bar x)t> 0 \; \forall  t \in  \Omega(W_{m_*})$ for   some  $\bar x\in X$.
\end{itemize}

The algorithm is described.$  \qquad \blacksquare$

\begin{remark} In the described above  algorithm \textbf{REG-LCoP}, it is assumed  that $X\not =\emptyset $.   It is easy to modify  the algorithm so that this assumption  is removed.
\end{remark}

\subsubsection{ On the  comparison   of the algorithms  }\label{subsec7}

To give another interpretation of the algorithm \textbf{REG-LCoP} and   to better trace the compliance of  the algorithm \textbf{REG-LCoP} to   the  Waki and  Muramatsu's FRA  from \cite{Waki} (presented in \ref{subsec5}), let us perform some additional constructions   at  the  iterations of the algorithm \textbf{REG-LCoP}.

\vspace{3mm}

At the end of
{\bf Iteration $\#$ 0}, having data $\tau(i),\;  \gamma(i), \; L_1(i), i\in I_1,$
 let us set
$${\mathcal F}_0{:}={\cal COP}^p,\;\; Y_1{:}=\sum\limits_{i\in I_1}\gamma(i)\tau(i)
(\tau(i))^{\top},\; {\mathcal F}_1{:}={\mathcal F}_0\cap \{Y_1\}^\bot.$$
Notice here that, by construction,
\begin{equation*}\begin{split}&\qquad\qquad\qquad Y_1\not =\mathbb O_{p},\; Y_1\in {\rm Ker }\mathcal  A,\; Y_1\in {\mathcal F}^*_0={\cal CP}^p,\\
{\mathcal F}_1: &={\mathcal F}_0\cap \{Y_1\}^\bot=\{D\in {\cal COP}^p: D\bullet Y_1=0\}=\{D\in {\cal COP}^p:(\tau(i))^\top D\tau(i)=0\  \forall  i \in I_1\}\qquad
\\&=\{D\in {\cal COP}^p:e^\top _kD\tau(i)=0 \  \forall  i \in {L_1(i)},\  e^\top _kD\tau(i)\geq 0\  \forall  i \in P\setminus {L_1(i)},\  \forall  i \in I_1\},\end{split}\end{equation*}
where $\mathbb O_p$ is the $p\times p$  null matrix.

\vspace{3mm}

Consider  {\bf Iteration $\#$ $m$, $1\leq m\leq m_*$.}
By the beginning of  the iteration, we have 
a cone $ {\mathcal F}_m={\mathcal F}_{m-1}\cap \{Y_m\}^\bot$ that  can be described as follows:
\begin{equation}\begin{split}{\mathcal F}_m
=\{D\in {\cal COP}^p: \ e^\top _kD\tau(i)=0 \  \forall  k \in L_m(i), \; e^\top _kD\tau(i)\geq 0\   \forall  k \in P\setminus L_m(i), \forall  i \in I_m\}.\label{o-14}\end{split}\end{equation}

At the end of this iteration, we   have new data (\ref{o-18})   and numbers $\gamma(i), i\in \Delta I_m$. Let us  set
\be Y_{m+1}{:}=\sum\limits_{i\in \Delta I_m}\gamma(i)\tau(i)(\tau(i))^{\top}+
\sum\limits_{i\in I_m} [\tau(i)(\lambda^m(i))^{\top}+\lambda^m(i)(\tau(i))^{\top}].\label{Ym}\ee
From the equations in (\ref{18-11}), we conclude:  $Y_{m+1}\in {\rm Ker}\mathcal  A.$  From (\ref{o-14}), it follows:
\bea &{\mathcal F}^*_m =\ {\rm cl}\{D\in \mathcal S(p):\ D\in {\cal CP}^p \oplus  {\cal P}^*_m \},\nonumber\\&
 {\cal P}^*_m:= \{D\in \mathcal S(p):D=\sum\limits_{i\in I_m}
[\tau(i)(\lambda(i))^{\top}+\lambda(i)(\tau(i))^{\top}],\ \lambda_k(i)\geq 0 \  \forall  k \in P\setminus {L_m(i)},  \forall  i \in I_m\}.\nonumber\eea
Hence, by construction, $Y_{m+1}\in {\mathcal F}^*_m.$

Consider the cone ${\mathcal F}_{m+1}:={\mathcal F}_{m}\cap \{Y_{m+1}\}^\bot$ and show that it
 can be described as follows{:}
\begin{equation}\begin{split}  {\mathcal F}_{m+1}
=\{D\in {\cal COP}^p:\  &e^\top _kD\tau(i)=0  \  \forall  k \in {L_{m+1}(i)},\\ &  e^\top _kD\tau(i)\geq 0 \   \forall  k \in P\setminus {L_{m+1}(i)}, \forall  i \in I_{m+1}\}.\label{o-15}\end{split}\end{equation}
In fact, it follows from (\ref{Ym}) that the  equality $D\bullet Y_{m+1}=0$ can be rewritten in the form
\begin{eqnarray} &0=D\bullet Y_{m+1}=\sum\limits_{i\in \Delta I_m}\gamma(i)(\tau(i))^{\top}D\tau(i)+2\sum\limits_{i\in I_m}(\lambda^m(i))^{\top}D\tau(i), \nonumber\\
&\mbox{ where }  \ \  \gamma(i)>0 \  \forall  i \in \Delta I_m;\;\;\;  \lambda^m_k(i)\geq 0\  \forall  k \in P\setminus {L_m(i)},  \forall   i \in I_m.\nonumber\end{eqnarray}
Taking into account (\ref{o-14}) and {the} relations above, we conclude that
for $D\in {\mathcal F}_m$, the  equality $D\bullet Y_{m+1}=0 $ implies the equalities
$$ (\tau(i))^{\top}D\tau(i)=0 \  \forall i \in \Delta I_m,\; e^\top_kD\tau(i)=0\   \forall  k \in \Delta {L}(i),  \forall  i \in I_m.$$
Notice that the relations $D\in {\cal COP}^{p},$ $(\tau(i))^{\top}D\tau(i)=0,  \tau(i)\geq 0  \ \forall  i \in \Delta I_m,$ imply
\be  e^\top_kD\tau(i)=0 \  \forall  k \in P_+(\tau(i)){\mbox{ and } }  e^\top_kD\tau(i)\geq 0\  \forall  k \in P\setminus  P_+(\tau(i)), \forall   i \in \Delta I_m.\label{o-17}\ee
Representation (\ref{o-15}) follows from (\ref{o-14})
  and (\ref{o-17}).

The  constructed by rules (\ref{Ym}) and (\ref{o-15})  matrices and cones
\be Y_m,\; {\mathcal F}_m,\; m=0, 1,...,m_*,\label{w4}\ee
satisfy the following relations:
\begin{equation*}\begin{split}&Y_m\in {\mathcal F}^*_{m-1}\  \ \  \forall  m=1,...,m_*, \ Y_0=\mathbb O_{p}; \qquad\qquad\qquad\qquad\qquad\qquad(I)\\
&Y_m\in {\rm Ker}\mathcal  A\ \ \   \forall  m=1,...,m_*; \qquad\qquad\qquad\qquad\qquad\qquad\qquad\qquad(II)\\
&{\mathcal F}_m= {\mathcal F}_{m-1}\cap \{Y_m\}^\bot\ \ \   \forall  m=1,...,m_*, \ {\mathcal F}_0={\mathcal {COP}}^p.\qquad\qquad\qquad(III)\end{split}\end{equation*}

Now we see that
 the algorithm {\textbf{REG-LCoP}}  allows one   to get   a  more clear description of the structure of the matrices $Y_m, m=1,...,m_*,$ satisfying conditions $(I)-(III)$,  and quite constructive rules of their formation:

  \begin{itemize}
  \item[]\it for a given $m$, the matrix $Y_m$  has a form  (\ref{Ym}) and is built on the  basis of {the} optimality conditions for the feasible solution  $(x=\mathbf 0, \mu=0)$ in the {corresponding} \textbf{regular}  SIP  problem (${SIP}_m$).
  \end{itemize}

 As it was shown in  subsection \ref{subsec5}, at each iteration,   the Waki and  Muramatsu's FRA produces  a set of matrices and cones (\ref{w4}) satisfying the {conditions}  $(I) - (III)$, and  the condition
$$Y_m\not\in {\rm span}\{Y_0,Y_1,...,Y_{m-1}\}\   \forall    m=1,...,m_*-1.\qquad\qquad\qquad (IV)$$

 On the other hand, t he described in subsection \ref{subsec6} algorithm \textbf{REG-LCoP}  at each iteration produces  a set of matrices and cones (\ref{w4})  satisfying the conditions  $(I) - (III)$  but not necessarily the condition (IV).

Since in the algorithm \textbf{REG-LCoP}, the fulfillment of   the condition (IV)  is  not guaranteed   at  each iteration,
if compare this algorithm  with the  Waki and Muramatsu's FRA, at the first glance it may seem that, in general,  the number of iterations executed by the algorithm \textbf{REG-LCoP} is larger.  Such an impression is caused by the fact that  in \ref{4-2-1},  we described in  more   detail all the steps of the algorithm and explicitly indicated all the computations  carried out at each iteration.  As for the  Waki and  Muramatsu's FRA, its  iterations  are described only in general terms.

In what follows, we   set out  a {\it modification} of the algorithm \textbf{REG-LCoP}, where  the number of  iterations  is reduced   and it is  guaranteed that all conditions $(I) - (IV) $ are satisfied on each {\it core} iteration.
This modification is formal, being essentially  another way of the iterations' numbering. The real number of the calculations on the steps of this modified algorithm  is the same as on the iterations of the original  one.

\subsubsection{ A compressed modification  of the algorithm \textbf{REG-LCoP} }

 Consider the algorithm \textbf{REG-LCoP}   presented  in  \ref{4-2-1}. Evidently,  one can reduce the number of  iterations of the algorithm if squeeze  into just one iteration   that iterations  of the algorithm which  change the description of the  dual cone  ${\cal F}^*_m$,
but do not   change the cone ${\cal F}_m$ itself. In other words, we will only move to the next {\it core} iteration  when all  conditions $(I) - (IV)$ are satisfied.  Formally,  such a procedure can be   described as follows.

\vspace{3mm}

Suppose that the algorithm \textbf{REG-LCoP} has constructed matrices and cones (\ref{w4}),  satisfying the properties $ (I)-(III)$   and let $m_*>0.$
Denote by
$m_s\in \{0,1,...,m_*-1\},\; s=0,1,...,s_*, $
the iterations'  numbers such that
\begin{equation*}\begin{split}&m_0: =0,\; m_s<m_{s+1}\;  \forall  s=0,1,...,s_*-1;\;\; m_{s_*}+1=m_*,\\
&Y_{m_s+1}\not\in {\rm span}\{Y_0,Y_1,...,Y_{m_s}\}\;  \forall  s=0,1, ...,s_*-1;\\
&Y_{m_s+1+i}\in {\rm span}\{Y_0,Y_1,...,Y_{m_s+i}\}\;  \forall  i=1,...,m_{s+1}-m_s-1, \ \forall  s=0,1,...,s_*-1.\end{split}\end{equation*}
Here $s_*$ denotes the number of  iterations  for which  the conditions above  are met.
Notice that the set $\{l,l+1,...,w\}$ is considered empty if $w<l.$

\vspace{2mm}

In other words, the condition $(IV)$ {is} satisfied only for
$m\in \{m_s+1,\; s=0,1,...,s_*-1\}$  and, possibly, for $m=m_*$. Set
$$\bar Y_0:=Y_0,\; \bar {\mathcal F}_0:={\mathcal F}_0,\; \bar Y_{s+1}:=Y_{m_s+1},\; \bar {\mathcal F}_{s+1}:={\mathcal F}_{m_s+1}\;  \forall  s=0,1,...,s_*.$$
It is easy to check that the following conditions hold true:
$$\bar Y_s\in \bar {\mathcal F}^*_{s-1},\; \bar Y_s\in {\rm Ker}\mathcal A,\; \bar {\mathcal F}_s=\bar {\mathcal F}_{s-1}\cap \{\bar Y_s\}^\bot \  \forall  s=1,...,s_*+1;$$
$$\bar Y_s\not\in{\rm span}\{\bar Y_0,\bar Y_1,...,\bar Y_{s-1}\}\;  \forall  s=1,...,s_*.$$

Thus, after the described above squeezing, we get $s_*$ core iterations of the modified algorithm.
It follows from the conditions above that $s_*\leq {\rm dim}(\rm Ker \,\mathcal A).$

\vspace{2mm}

Notice that
for any  $s=0,1,...,s_*-1,$  the  iterations of the algorithm \textbf{REG-LCoP} having the numbers
$m_s+1+i,\; \mbox{where }  i=1,...,m_{s+1}-m_s-1$ (the compressed iterations),  are not
useless. They
can be considered as the steps of a regularization
procedure for the cone ${\mathcal F}_{m_s+1}$ at the current  core iteration $\# \; s $.
At each of these iterations, we reformulate   the cone ${\mathcal F}_{m_s+1}$  in a new equivalent form. This  additional information allows us to improve (make more regular) the  representation of  the cone ${\mathcal F}_{m_s+1}$ and get a more explicit and  useful  description of its  dual cone  ${\mathcal F}^*_{m_s+1}$.

\vspace{2mm}

\subsubsection{A short discussion on the algorithms considered in this section}\label{Discussion}

By analyzing and comparing the iterative algorithms presented above, we can draw the following conclusions.

\begin{enumerate}
\item  The  Waki and  Muramatsu's  facial reduction algorithm from \cite{Waki}, reformulated for copositive problems in  subsection \ref{subsec5}, is  very simple in  the  description  and runs no more  than ${\rm dim } ({\rm Ker}\mathcal A)$ iterations. But this algorithm is more   conceptual  than constructive since it does not provide  any information about the structure of the matrix  $Y_m$ and  the cone  ${\mathcal F}^*_m$ at its $m$-th iteration. Moreover, it is not explained in \cite{Waki} how to fulfill  steps 2 and 3 at each iteration.

\item The algorithm \textbf{REG-LCoP} proposed in  subsection \ref{subsec6} also runs a finite number of iterations.   This algorithm is    described   in all details and justified.  The  quite constructive rules for calculating the matrix $Y_m$ satisfying the condition   $Y_m\in {\mathcal F}^*_m$, are presented using  the  information available at the Iteration $\# m$ of this algorithm. These rules are derived from  the optimality conditions  for  the  optimal solution  $(x=\mathbf 0, \mu=0)$ of the \textbf{regular} problem (${SIP}_{m})$.

Notice that it is possible to develop a modification of the algorithm  \textbf{REG-LCoP}  which runs no more than $2n$ iterations.

\item Finally,  to show that the  described in  \ref{4-2-1}  algorithm  \textbf{REG-LCoP}   is not worse (by the number of iterations) than the FRA from the subsection \ref{subsec5}, we  presented  a  compressed modification of the algorithm  \textbf{REG-LCoP}. This modification  consists of no more than ${\rm dim } ({\rm Ker}\mathcal A)$ iterations  as well as the algorithm from subsection \ref{subsec5}.

\end{enumerate}

\section{Conclusions}\label{Concl}

The main contribution  of the paper    is that, based on the concept of immobile indices, previously introduced for semi-infinite  optimization problems, we  suggested new  methods for  regularization of  copositive problems. The algorithmic  procedure  of  regularization  of   copositive problems is described in the form of  the algorithm \textbf{REG-LCoP} and is compared  with the  facial reduction approach  based on the minimal cone representation.  We show that, when applied to the linear CoP problem (\ref{SDP-2}), the  algorithm \textbf{REG-LCoP}  possesses    the same properties as the FRA  suggested by  Waki and Muramatsu in \cite{Waki}, but its iterations are explicit,  described in more detail and hence more constructive.

The   described  in  the paper algorithms are useful for  the  study of convex copositive problems. In particular, for  the  linear copositive problem, they allow  to
    \begin{itemize}
\item    formulate an    equivalent (regular) semi-infinite problem which satisfies the Slater type regularity condition and can be solved numerically;
    
    \vspace{-2mm}
    
     \item    prove   new optimality conditions without  any CQs;
     
     \vspace{-2mm}
     
     \item    develop  strong duality theory  based on an  explicit representation of the  "regularized"  feasible cone and the corresponding dual (such as, e.g. the Extended Lagrange Dual Problem suggested for SDP by  Ramana et al. in \cite{Ramana-W}).
\end{itemize}
The described in the paper regularization approach is novel and rather constructive. It is important to stress that  no  constructive regularization procedures are known for linear copositive problems.

 \section*{Acknowledgement} This work was partially supported  by the  state research program "Convergence"(Republic Belarus), Task 1.3.01, by
 Portuguese funds through CIDMA - Center for Research and Development in Mathematics and Applications, and FCT -
 Portuguese Foundation for Science and Technology, within the project  UIDB/04106/2020.

\end{document}